\documentclass[12pt]{amsart}
\usepackage{amsmath}
\usepackage{amsthm}
\usepackage{amscd}

\newcommand{\Q}{{\mathbb Q}}

\newcommand{\Z}{{\mathbb Z}}
\newcommand{\C}{{\mathbb C}}
\newcommand{\F}{{\mathbb F}}
\newcommand{\A}{{\mathbb A}}

\newcommand{\fg}{{\mathfrak g}}

\newtheorem{theorem}{Theorem}%[section] (If you want theorem numbered
\newtheorem{lemma}{Lemma}%               with section number.  Same
\newtheorem{corollary}{Corollary}%       goes for lemmas, etc.)
\newtheorem{proposition}[theorem]{Proposition}
\numberwithin{theorem}{section}
\numberwithin{corollary}{section}
\numberwithin{lemma}{section}

\theoremstyle{definition}

\numberwithin{conj}{section}
\newtheorem{example}{Example}
\newtheorem*{acknowledgement}{Acknowledgement}
\numberwithin{example}{section}
\newtheorem{definition}{Definition}
\numberwithin{definition}{section}
\newtheorem{question}{Question}
\numberwithin{question}{section}
\numberwithin{equation}{section}

\theoremstyle{remark}
\newtheorem{remark}{Remark}
\numberwithin{remark}{section}
\newtheorem*{notation}{Notation}

\begin{document}

\title{Recovering $l$-adic representations}
%\date{30th November, 2001}
\author{C.~S.~Rajan}

\address{Tata Institute of Fundamental 
Research, Homi Bhabha Road, Bombay - 400 005, INDIA.}
\email{rajan@math.tifr.res.in}


\begin{abstract}
We  consider  the problem of 
recovering  $l$-adic representations from a
knowledge of the character values at the Frobenius elements associated
to $l$-adic representations constructed algebraically out of the original
representations. These results generalize earlier results 
in  \cite{rajan98} concerning refinements of strong multiplicity one for
$l$-adic represntations,  and a result of Ramakrishnan
\cite{DR} recovering modular forms from a knowledge of the
squares of the Hecke eigenvalues.    
For example, we show that if the characters of some  tensor or
symmetric powers of two absolutely irreducible $l$-adic representation
with the algebraic envelope of the image being connected, agree 
at the Frobenius elements corresponding to a set of places of
positive upper density, then the representations are twists of each
other by a finite order character.
\end{abstract}

\maketitle

\section{Introduction}
Let $K$ be a global field and let $G_K$ denote the Galois group over
$K$ of an algebraic closure $\bar{K}$ of $K$. Let $F$ be a
non-archimedean local field of characteristic zero, and $M$ an affine,
algebraic group over $F$. Suppose 
 $$\rho_i:G_K\rightarrow M(F), ~~i=1,2$$
are continuous, semisimple  representations 
of the Galois group $G_K$ into $GL_n(F)$, unramified outside a finite
set $S$ of places containing the archimedean places of $K$.    
 Let 
\[ R: ~M\to GL_m\]
be a rational  representation of $M$ into $GL_m$ defined over $F$. We
assume that the kernel $Z_R$ of $R$, is contained inside the centre
of $M$.
In this paper,  we are interested in the problem  of recovering $l$-adic
representations from a knowledge of $l$-adic representations
constructed algebraically out of the original representations, as follows:

\begin{question}\label{question} Let $T$  be a  subset  of the set of 
places $\Sigma_K$ of $K$ satisfying,      
\begin{equation}\label{T}
 T=\{v \not\in S \mid{\rm Tr}({R\circ\rho_1}(\sigma_v))= {\rm Tr}({R\circ
\rho_2}(\sigma_v))\}, 
\end{equation}
where $\rho_i(\sigma_v)$ for $v\not\in S$
denotes the Frobenius conjugacy classes lying in
the image. Suppose that $T$ is a  `sufficiently large' set of places of $K$. 
  How are $\rho_1$
and $\rho_2$ related? 

More specifically, under what conditions on $T$,  $R$ or the nature
of the representations $\rho_i$, can we conclude that there exists a
central abelian representation $\chi: G_K\to Z_R(F)$, such that the
representations $\rho_2$ and $\rho_1\otimes \chi$ are conjugate by an
element of $M(F)$?

Further  we would also like to know the answer when we take $R$ to be a
`standard representation', for example if $R$ is taken to be $k^{th}$
tensor, or symmetric or exterior powers of a linear representation of
$M$.  The representation $R\circ\rho$ can be thought of 
as an $l$-adic representation
constructed algebraically from the original representation $\rho$.  

\end{question}

When $M$ is isomorphic to $GL_m$,
 and $R$ is taken to be the identity morphism, then the
question is a refinement of strong multiplicity one and was considered
in  earlier papers \cite {DR1}, \cite{rajan98}. In \cite{rajan98}, 
we proved Ramakrishnan's
conjecture \cite{DR1}, that if the upper density of
$T$ is strictly greater than $1-1/2m^2$, then $\rho_1\simeq
\rho_2$. Further,  if the upper density of $T$ is
sufficiently large in relation with the number of connected components
of the algebraic envelope of the image of $\rho_1$, then 
 the representations are isomorphic upon  restriction  to the absolute
Galois group of a  finite extension of $K$. Morever it was shown that    
if $\rho_1$ is absolutely irreducible, the algebraic envelope of the
 image of 
$\rho_1$ is connected and $T$ is of positive density, then there
exists a character $\chi$ of $G_K$, such that $\rho_2\simeq
\rho_1\otimes \chi$ (Theorem \ref{smo} below).

 It was shown by Ramakrishnan
\cite{DR}, that holomorphic newforms are determined upto a quadratic
twist from knowing that the squares
of the Hecke eigenvalues coincide at a set of places of density at
least $17/18$.  In this paper we use  the 
methods and techniques of our earlier
paper \cite{rajan98},  to obtain generalisations in the
context of $l$-adic representations, of the theorem of
Ramakrishnan.  In the above setting, 
Ramakrishnan's theorem can be obtained by taking $M$ to
be $GL_2$,  and $R$ to be either the symmetric square or the two
fold tensor product of the two dimensional regular representation of
$GL(2)$. In the process we obtain a different proof (and also a
generalisation)  of
Ramakrishnan's theorem.

We  break the general problem outlined above  into two steps.  First we use the
results of \cite{rajan98}, to conclude that $R\circ \rho_1$ and 
$R\circ \rho_2$ are isomorphic under suitable density hypothesis on $T$.
We consider then the algebraic envelopes of the $l$-adic
representations, and  we try to answer the question of recovering
algebraic representations  with the role
of the Galois group replaced by a  reductive group. 
We first state a general theorem that under suitable density hypothesis on $T$, we can  conclude upto
twisting by a central abelian representation, that the representation $\rho_2$ is essentially
an `algebraic twist' of  $\rho_1$. Combining this with a theorem of
Richardson \cite{Ri}, we obtain a finiteness result, that if we fix a
$l$-adic representation $\rho_1$,  then  the collection of
$l$-adic representations $\rho_2$ satisfying the hypothesis of the
above question for some $R$ 
fall into finitely many `algebraic classes' upto twisting by an
abelian representation commuting with $\rho_1$. 

We then specialize $R$ to be either  symmetric, tensor power, adjoint and
twisted product (Asai) 
representations of the ambient group $GL_n$. In this situation, the
ambiguity regarding algebraic twists can be
resolved, and under suitable density hypothesis, we can essentially recover 
representations from knowing that the character values of the
 symmetric or tensor powers of them coincide.  We discuss the
situation when the representation is absolutely irreducible and the
algebraic envelope of the image is not connected, with the expectation
that  it may be of use in understanding some aspects of endoscopy, as
for instance arising in the work of Blasius \cite{blasius94}. 

Specialising to  modular forms we generalise the results of
 Ramakrishnan \cite{DR},   where we
consider arbitrary $k^{th}$ powers of the eigenvalues for a natural
number $k$, and also $k^{th}$ symmetric powers in the eigenvalues of
the modular forms, and the Asai representations. 
 We also present an  application  in the context of  abelian
varieties. 

\begin{acknowledgement} The initial ideas for this paper was conceived
when the author was visiting Centre de Recherches Math\'{e}matiques,
Montr\'{e}al during the Special Year on Arithmetic held in 1998, and
my sincere  thanks to  CRM for their  hospitality and support. I thank
M. Ram Murty for useful discussions and 
the  invitation to visit CRM, and to Gopal Prasad for
helping me  with the references to the work of Richardson's and  of Fong and
Greiss. Some of these results were indicated in \cite{Ra2}.
 The arithmetical 
application of  Theorem \ref{tensorandsym} d),
to generalized Asai representations was suggested by D. Ramakrishnan's
work, who had earlier proved a similar result for the usual degree two
Asai representations, and I thank him for conveying to me his
results.

\end{acknowledgement}

\section{Preliminaries}

\begin{notation} 
Denote by  $\Sigma_K$ the set of places of $K$. 
For a nonarchimedean place $v$ of $K$,  let ${\mathfrak p}_v$ denote 
the corresponding prime ideal of ${\mathcal O}_K$,  and
$Nv$ the norm 
of $v$ be the number of elements of the finite field 
${\mathcal O}_K/{\mathfrak p}_v$. Suppose  $L$ is  a finite Galois extension of $K$,
with Galois group $G(L/K)$. Let $S$ denote a  subset of $\Sigma_K$, 
containing the archimedean places together with the set of places of $K$
which ramify in $L$. For each place $v$ of $K$ not in $S$, and a place 
$w$ of $L$ lying over $v$, we have a canonical Frobenius element $\sigma_w$
in $G(L/K)$, defined by the following property:
\[ \sigma_w(x)\cong x^{Nv}({\rm mod}\,{\mathfrak p}_w).\]
 The set $\{\sigma_w\mid w|v\}$ form the Frobenius conjugacy 
class in $G(L/K)$, which we continue to denote by $\sigma_v$. 

  Let $M$ denote a connected, reductive 
algebraic group defined over $F$, and let $\rho$ be
 a continous representation of $G_K$ into 
$M(F)$, where $F$ is a non-archimedean 
local field of residue characteristic $l$. Let $L$ denote the fixed field
of $\bar{K}$ by the kernel of $\rho.$ Write $L=\cup_{\alpha}L_{\alpha},$ where 
$L_{\alpha}$ are finite extensions of $K$. $\rho$ is said to be unramified 
outside a set of primes $S$ of $K$, if each of the extensions 
$L_{\alpha}$ is an unramified extension of $K$ outside $S$.

We will assume henceforth that all our linear $l$-adic 
representations of $G_K$  are continuous and semisimple, since we need
to determine a linear representation from it's character.  By the
results of \cite{KR}, 
 it follows that the set of ramified primes is of
density zero, and hence arguments involving density as in
\cite{rajan98}, go through essentially unchanged. $S$ will indicate a
a set of primes of density zero, containing the ramified primes of the
(finite) number of l-adic representations under consideration, and the
archimedean places of $K$. 
 
 Let $w$ 
be a valuation on $L$ extending a valuation $v\not\in S.$ The Frobenius 
elements at the various finite layers for the valuation $w\mid_{L_{\alpha}}$ 
patch together to give raise to the Frobenius element $\sigma_w\in G(L/K)$,
and a Frobenius conjugacy class $\sigma_v\in G(L/K).$ Thus $\rho(\sigma_w)$
(resp. $\rho(\sigma_v)$) is a  well defined element (resp. 
conjugacy class) in $M(F)$. 
If $\rho:G_K\to GL_m(F)$ is a linear representation,
let $\chi_{\rho}$ denote the associated character. 
 $\chi_{\rho}({\sigma_v})$ is well
defined for $v$ a place of $K$ not in $S$.

 Let $G$ be a connected, reductive group. We have a decomposition,
\[ G=C_G G_d,\]
where $C_G$ is the connected component of the center of $G$,  $G_d$ the
semisimple component of $G$, is
the derived group of $G$,  and
$C_G\cap G_d$ is a finite central subgroup of $G_d$.  

For a diagonalisable group  $D$, let $X^*(D)$ denote the finitely
generated dual abelian  group
of characters of $D$. $D$ is a torus if and only if $X^*(D)$
is a free abelian group. Given a morphism $\phi:D\to D'$ between 
diagonlisable groups, let $\phi^*$ denote the dual morphism
$X^*(D')\to X^*(D)$. We recall that there is an (anti)-equivalence
between the category of diagonalisable groups over $\bar{F}$ 
and the category of
finitely generated abelian groups, where the morphism one way 
at the level of objects is given by sending  $D$ to it's character
group $X^*(D)$. 

 For a continuous morphism
 $\rho: G_K\to M(F)$, denote by  $H_{\rho}$ 
denote the algebraic envelope in $M$ of  the image group $\rho(G_K)$ inside
 $M(F)$,  i.e., the 
smallest algebraic subgroup $H_{\rho}$ of $M$, defined over $F$, such that 
$\rho(G_K)\subset H_{\rho}(F)$. $H_{\rho}$ is also the Zariski closure
over $F$ of $\rho(G_K)$
 inside $M$. 
For $i=1,2$, let $H_i=H_{\rho_i}$, and let $H_i^0$ be the identity
 component of $H_i$. 

We let $R:M\to GL_m$ denote a rational representation of $M$, and $Z$
a central subgroup of $M$ containing the kernel $Z_R$ of $R$. 
\end{notation}

\section{Strong multiplicity one}
We recall here the results of \cite{rajan98}. Let $R:M\to GL_m$ be a
rational representation of $M$. The upper density  
$ud(P)$ of a set $P$ of primes of $K$, is defined to be the ratio,
\[
ud(P)={\varlimsup}_{x\rightarrow \infty}\frac
{\#\{v\in \Sigma_K\mid Nv\leq x,~v\in P\}}
{\#\{v\in \Sigma_K\mid Nv\leq x\}}.\]
Consider the following two hypothesis on the upper density of $T$,
depending on the representations $\rho_1$ and $\rho_2$:
\begin{align}
& {\bf DH1:} &  ud(T) &  >1-1/2m^2 \\
& {\bf DH2:} & ud(T) &  >{\rm min}(1-1/c_1,1-1/c_2),
\end{align}
where $c_i=|R(H_i)/R(H_i)^0|$
 is  the number of connected components of $R(H_i)$.
As a consequence of the refinements of strong multiplicity one proved
in \cite[Theorems 1 and 2]{rajan98}, we obtain
\begin{theorem} \label{smo} 
 i) If $T$ satisfies {\bf DH1},  then $R\circ\rho_1\simeq R\circ\rho_2$.

ii) If $T$ satisfies {\bf DH2},
then  there is a finite Galois extension $L$ of $K$, such that 
$R\circ\rho_1\!\mid_{G_L}\simeq R\circ\rho_2\!\mid_{G_L}$.
 The connected component $R(H_2^0)$  is conjugate to $R(H_1^0)$. In
particular if 
either $H_1$ or $H_2$ is connected and $ud(T)$ is positive, then there
is a finite Galois extension $L$ of $K$, such that 
$R\circ \rho_1\!\mid_{G_L}\simeq R\circ \rho_2\!\mid_{G_L}$.

\end{theorem}
\begin{proof}
For the sake of completeness of exposition, we present a brief outline
of the proof and refer to \cite{rajan98} for more details. The proof
reduces to taking $M=GL_m$ and $R$ to be the identity
morphism. Let $\rho=\rho_1\times \rho_2$, and let $G$ denote the image
$\rho(G_K)$. Consider the algebraic subscheme 
\[ X=\{(g_1,g_2)\mid {\rm Tr}(g_1)={\rm Tr}(g_2)\}.\]
It is known that if $C$ is a closed, analytic subset of
$G$, stable under conjugation by $G$ and of dimension strictly smaller
than that of $G$, then the set of Frobenius conjugacy classes lying in
$C$ is of density $0$. Using this it follows that the collection of 
Frobenius conjugacy classes lying in
$X$ has a density equal to 
\begin{equation}\label{algdensity}
\lambda= \frac{|\{\phi\in H/H^0\mid H^{\phi}\subset X\}|}{|H/H^0|}.
\end{equation}
Since this last condition is algebraically defined, the above
expression can be calculated after base changing to $\C$. Let $J$
denote a maximal compact subgroup of $H(\C)$, and let $p_1,~p_2$
denote the two natural projections of the product $GL_m\times GL_m$. 
Assume that $p_1$ and $p_2$ give raise to inequivalent representations
of $J$. (i) follows from the inequlities
\begin{equation}\label{orthogonality}
 2\leq \int_J |{\rm Tr}(p_1(j))-{\rm Tr}(p_2(j))|^2d\mu(j)\leq
(1-\lambda)4m^2,
\end{equation}
where $d\mu(j)$ denotes a normalized Haar measure on $J$.
 The first
inequality follows from the orthogonality relations for
 characters. For the second inequality, we observe that the eigenvalues of
 $p_1(j)$ and $p_2(j)$ are roots of
 unity, and hence 
\[ |{\rm Tr}(p_1(j))-{\rm Tr}(p_2(j))|^2\leq 4m^2.\]
Combining this with the expression for the density $\lambda$ given by
equation \ref{algdensity} gives us the second inequality.

To prove (ii), it is enough to show that $H^0\subset X$. Let
$c_1<c_2$. The density hypothesis implies together with the expression
(\ref{algdensity}) for the density, that there is some element of the
form $(1,j)\in J\cap X$. The proof concludes by observing that the
only element in the unitary group $U(m)\subset GL(m,\C)$ with trace
equal to $m$ is the identity matrix, and hence the connected component
of the idenity in $J$ (or $H$) is contained inside $X$. 
\end{proof} 

\begin{remark}
In the automorphic context, assuming the Ramanujan-Petersson
conjectures, it is possible to obtain the inequalities in
(\ref{orthogonality}), by analogous arguments, and thus a proof of
Ramakrishnan's conjecture in the automorphic context. The first
inequality follows from replacing the orthogonality relations for
characters of compact groups, by the Rankin-Selberg convolution of
$L$-functions, and amounts to studying the behvior 
at $s=1$ of the logarithm  of the function, 
\[L(s,`{|\pi_1-\pi_2|^2}'):= \frac{L(s,\pi_1\times
\tilde{\pi}_1)L(s,\pi_2\times \tilde{\pi}_2)}{L(s,\pi_1\times
\tilde{\pi}_2)L(s,\tilde{\pi}_1\times {\pi}_2)},\]
 where $\pi_1$ and $\pi_2$ are unitary,
automorphic representations of $GL_n({\bf A}_K)$ of a number field $K$, and
$\tilde{\pi}_1, ~\tilde{\pi}_2$ are the contragredient representations
of $\pi_1$ and $\pi_2$. The second
inequality follows from the Ramanujan hypothesis. For more details we
refer to \cite{Ra1}. 
 \end{remark}

\section{Recovering representations}
 Theorem \ref{smo}  allowed  us to conclude under
certain density hypothesis on $T$, that  $R\circ \rho_1\simeq R\circ
 \rho_2$. Our aim now is to find sufficient algebraic and group
 theoretic conditions in order to deduce that $\rho_1$ and $\rho_2$
determine each other upto twisting by a central abelian
representation. The problem essentially boils down to considering
extensions of morphisms between tori.  

Let $\phi: H_1 \to H_2$ be an algebraic  homomorphism between 
  reductive groups defined over $F$. Suppose $H_2'$ is a 
reductive group over $F$, with a {\em surjective} homomorphism,
\[\pi:~H_2'\to H_2,\]
with $Z_{\pi}'={\rm Ker}(\pi)$ a central subgroup in  $H_2'$.    
We  say that a
homomorphism $\psi: H_1\to H_2'$ defined over $\bar{F}$, lifts $\phi$, if 
\[\pi\circ \psi =\phi.\] 
We have the following lemma giving sufficient conditions for a lift to
exist:

\begin{lemma}\label{lift} Suppose $\phi: H_1 \to H_2$ is a homorphism between 
connected,  reductive groups over $F$.
With notation as above, we have the following: 
\begin{enumerate}
\item There exists a connected, reductive group $H_1'$ together 
with a isogeny $f: H_1'\to H_1$, and a lift  $\psi:H_1' \to H_2'$
of $\phi\circ f:~H_1'\to H_2$.

\item Suppose the derived group $H_{1d}$ is
simply connected, and that the kernel $Z_{\pi}$ of $\pi$ is connected.  
Then  a lift $\psi$ of $\phi$ exists.  

In particular if the derived group $H_{1d}$ is
simply connected, $H_2=PGL(n), ~H_2'=GL(n)$, then a lift exists.  
\end{enumerate}
\end{lemma}

\begin{proof} 1)  Let $H_{1s}$ be the simply connected cover of the
derived group 
$H_{1d}$ of $H_1$.  Let $H_1''=C_{H_1}\times H_{1s}$. 
We have a natural finite morphism
$f'':H_1''\to H_1$. Since $H_{1s}$ is simply connected, a lifting 
$\psi_s: H_{1s}\to H_2'$  of $\phi\circ f''\!\mid H_{1s}$ exists. 
We have reduced the proof to the following statement:
\begin{itemize}
\item $H_1=C_1\times H_{1s}$, $H_{1s}$ simply connected and $C_1$ a
torus. 
\item There are morphisms $\psi_s: H_{1s}\to H_2'$ and $\phi_c:C_1\to
 H_2$, such that the image of $\phi_c$ commutes with the image of
 $\phi_s=\pi\circ\psi_s$.  
\item we need to produce a torus $C_1'$ together with a isogeny
 $f_c: C_1'\to C_1$ and a lift $\psi_c:C_1'\to H_2'$ of $\phi_c\circ
 f_c$, such that it's image is contained in the commutant of  the
 image group $\psi_s(H_{1s})$. 
\end{itemize}
Let $T_2$ be the image torus
$\phi_c(C_1)$. Since the kernel of $\pi$ is central, we can find a
torus 
$T_2'$ contained inside $H_2'$ which surjects via $\pi$ onto $T_2$,
and commutes with the image $\psi_s(H_{1s})$.  We
have an injection at the dual level $X^*(T_2)\xrightarrow{\pi^*}
X^*(T_2')$, and let 
\[ M(\pi)=\{ \alpha\in X^*(T_2')\mid n\alpha \in \pi^*(X^*(T_2))
~{\rm for ~some} ~n\in \Z\},\]
be the torsion closure of $\pi^*(X^*(T_2))$ inside $X^*(T_2')$.
$M(\pi)$ is a direct summand of $ X^*(T_2')$.
Form the diagram, 
\[\begin{CD}
X^*(T_2) @>{(\phi\circ f)^*}>> X^*(C_{H_1})\\
@V{\pi^*}VV  @VVV\\
M(\pi) @>>> M_1
\end{CD}
\]
where $M_1$ is  the maximal torsion-free subgroup of the  group 
\[\frac{M(\pi)\oplus X^*(C_{H_1})}{{\rm Image}(\pi^*\oplus (\phi\circ f)^*)(X^*(T_2))}.\] 
Since  $\pi^*(X^*(T_2))$ is of finite index in $M(\pi)$, 
$M_1$ is a free abelian group of the same rank as $X^*(C_{H_1})$.
By the duality theory for tori, let $C_1'$ be a torus with character
 group $M$. $C_1'$ comes equipped with a finite map
 to $C_{H_1}$ corresponding to the dual map on the character groups
$X^*(C_{H_1})\to M_1$. Since $M(\pi)$ is a direct summand of $
X^*(T_2')$, there is a map from $C_1'\to T_2'$, and by construction
this commutes with the image $\psi_s(H_{1s})$.

2) Let $\phi_d$ denote a lift of $\phi|H_{1d}$ to
$H_2'$. In order to produce a lift of $\phi$, we need to extend
the morphism $\phi_d|C_{H_1}\cap H_{1d}: C_{H_1}\cap H_{1d} \to
H_2'$
to a map from $C_{H_1}$ to $H_2'$ commuting with the image
$\phi_d(H_{1d})$. 
 Let $T$ be the image torus
$\phi(C_{H_1})$. Let $T'$ 
be the inverse image of $T$ with respect to
$\pi$. Since kernel of  $\pi$ is central, $T'$ commutes with the image
$\phi_d(H_{1d})$. Hence, it is enough to extend the map
$\phi_d|C_{H_1}: C_{H_1}\cap H_{1d} \to
T'$ to a map from $C_{H_1}$ to $T'$. 
By the duality theory for tori, to produce a lift over $\bar{F}$ 
is equivalent to
producing a map from $X^*(T')$ to $X^*(C_{H_1})$ making the
following diagram commutative: 
\[\begin{CD}
X^*(C_{H_1}\cap H_{1d}) @<<<X^*(T')\\
@AAA  @AAA\\
X^*(C_{H_1}) @<<< X^*(T)
\end{CD}
\]
We have the exact sequence of tori,
\[ 1\to Z_{\pi}\to T'\to T\to 1,\]
and the corresponding dual sequence,
\[  1\to X^*(T)\to X^*(T')\to X^*(Z_{\pi})\to 1,\]
where the surjection follows from the fact that characters of closed
subgroups of diagonalisable groups extend to the ambient group. 
Since $Z_{\pi}$ is connected, we have a splitting $s:X^*(Z_{\pi})\to
X^*(T')$ and 
\[ X^*(T')=X^*(T)\oplus s(X^*(Z)).\]
  We are reduced then to producing a lift from $s(X^*(Z_{\pi}))$
to $X^*(C_{H_1})$ compatible with the projection to $X^*(C_{H_1}\cap
H_{1d})$. But since $Z_{\pi}$ is connected, $s(X^*(Z_{\pi}))$
is free, and by duality we obtain a lift as desired.
 
\end{proof} 

\begin{remark}\label{rationality}{\em (Rationality)}
 The question arises of imposing sufficient conditions
in order to ensure that $\theta$ is defined over $F$. From the proof
given above, the obstruction to the lift being defined over $F$ lies
essentially at the level of the morphism restricted to the center. 
 One sufficient
condition is to assume that $C_{H_1}$ and $Z$ are split tori over $F$,
and $T'\simeq T\times Z$ over $F$. For example, if
$M=GL_n$, the representation $\rho_1$ is absolutely irreducible
with $H_1$ connected and $H_{1d}$ simply connected, and $Z$ is the group of scalar matrices, then 
$\theta$ can be taken to be defined over $F$. 
\end{remark}

\begin{definition} Let $\rho_1:G_K\to M(F)$ be a $l$-adic
representation as above, and let $H_1$ be the algebraic envelope of
the image of $\rho_1$.  A $l$-adic representation $\rho_2$ is an
{\em (algebraic) conjugate} of $\rho_1$, if there exists an algebraic
homomorphism $\theta:H_1\to M$ defined over $\bar{F}$, such that 
\[\rho_2=\theta\circ \rho_1.\]
Let $Z$ be a subgroup of the center of $M$. 
$\rho_2$ is said to be a conjugate of $\rho_1$ upto twisting by
a central abelian representation with values in $Z$, if there exists a
homomorphism 
 $\chi:G_K\to Z(\bar{F})$ such that for any $\sigma\in G_K$,  we have
\[\rho_2(\sigma)=(\theta\circ \rho_1)(\sigma)\chi(\sigma).\]
\end{definition}
Let $R:M\to GL_m$ be a rational representation of $M$ such that the
kernel $Z_R$ of $R$ is a central subgroup of $M$. Combining Lemma
\ref{lift} and Theorem \ref{smo}, we obtain the following
general theorem giving sufficient conditions to recover $l$-adic
representations under suitable density hypothesis for the
corresponding characters: 

\begin{theorem} \label{recoverthm}
Let $\rho_i:G_K\to M(F), ~i=1,2$ be $l$-adic
representations as above. 
With the above notation, we have the following:

i) Suppose that  $R\circ\rho_1$ and $R\circ\rho_2$ satisfy {\bf
DH2}. Then the following hold:
\begin{itemize}
\item  there exists  a connected, reductive group $H_1'$ over
$\bar{F}$  with a finite
morphism $f: H_1'\to H_1^0$ and an algebraic homomorphism $\theta:
H_1'\to M$.

\item  a finite extension $L$ of $K$, and a
splitting $\rho_1': G_L\to H_1'(\bar{F})$ satisfying $f\circ
\rho_1'=\rho_1\!\mid G_L$. 

\item there exists a representation  $\chi:
G_L\to Z_R(\bar{F})$ such that, 
\[ \rho_2\!\mid G_L=(\theta \circ \rho_1')\otimes \chi.\]
\end{itemize}

ii)    Suppose that $H_1$ is connected, and 
the semisimple component $H^0_{1d}$  is simply connected.
Let $Z$ be a connected, central subgroup of $M$
containing the kernel $Z_R$ of $R$. 

a) Suppose further that  $R\circ\rho_1$ and $R\circ\rho_2$ satisfy {\bf DH2}.
Then there exists a finite extension $L$ of $K$, such that
$\rho_2|G_L$ is an algebraic twist of $\rho_1|G_L$ upto twisting by a
representation with values in $Z$.  

b)  Suppose that  $R\circ \rho_1$ and $R\circ \rho_2$ satisfy {\bf DH1}.
Then $\rho_2$ is an algebraic twist of $\rho_1$ upto twisting by a
representation with values in $Z$.

\end{theorem}

\begin{proof} By Theorem \ref{smo}, we can assume by the density
hypothesis {\bf DH2}, that there
exists a finite extension $L$ of $K$, such that 
\[ R\circ \rho_1|G_L\simeq R\circ \rho_2|G_L.\]
If {\bf DH1} is satisfied, we can take $L=K$. 

  For $i=1,2$, let
$\bar{H}_i$ denote the image of $H_i$ in the quotient group $M/Z_R$. The
morphism $R$ factors via 
$$R: M\xrightarrow{\pi} M/Z_R \xrightarrow{\bar{R}} GL_m,$$
where $\pi:M\to M/Z_R$ denotes the projection map, and $\bar{R}$ is the
embedding of $M/Z_R$ into $GL_m$ given by $R$.   Denote by 
\[\bar{\rho}_1, \bar{\rho}_2: G_K\to (M/Z_R)(F)\]
obtained from  the inclusion $M(F)/Z_R(F) \subset (M/Z_R)(F)$.  
Since $R\circ\rho_1\simeq R\circ\rho_2$, there is $A\in GL_m(F)$ such that 
\[ R\circ\rho_2(g)= A^{-1}(R\circ\rho_1)(g)A, ~~~~~~~~~~~~~\forall g\in G_K. \]
Denote by $\theta_A$ the automorphism of $GL_m$, given by conjugation
by $A$. The above equation translates to,
\[ \theta_A\circ\bar{R}\circ\bar{\rho}_1=\bar{R}\circ\bar{\rho}_2 \]
as morphisms from $G_K$ to $GL_m(F)$. The image
subgroups $\bar{\rho}_1(G_K)$ and $\bar{\rho}_2(G_K)$ are Zariski
dense in $\bar{H}_1$ and $\bar{H}_2$ respectively. Thus $\theta_A\circ
\bar{R}$ defines a homomorphism from $\bar{H}_1$ to
$\bar{R}(\bar{H}_2)$. Since $\bar{R}$ is an isomorphism of $M/Z_R$ onto
its image,  there exists an
algebraic isomorphism
\[ \bar{\theta}: \bar{H}_1\to \bar{H}_2, \]
satisfying \[ \bar{\theta}\circ\bar{\rho_1} =\bar{\rho_2}.\]

To prove the first part, by Lemma \ref{lift}, 
there exists a connected reductive group $H_1'$
over $\bar{F}$, a
finite morphism $f: H_1'\to H_1^0$,  and a morphism $\theta: H_1'\to
H_2$ lifting $\bar{\theta}\circ {\pi_1}$, where $\pi_1$ is the
projection from $H_1\to \bar{H}_1$.  Since $f$ is a finite central isogeny,
 pulling back this central extension of $H_1^0(F)$ to $G_K$ by the
morphism $\rho_1$, we obtain a finite, central extension of
$G_K$ by a finite abelian group $A$.  This gives an element $\alpha$ in
$H^2(G_K, A):={\rm lim}_{\xrightarrow{L}} H^2(G_K/G_L, A)$, where $L$
runs over the collection of finite Galois extensions of $K$ contained
inside $\bar{K}$. Thus there is a finite extension $L$ of $K$, at
which $\alpha$ can be considered as a  class in $H^2(G_K/G_L, A)$, and
then the restriction of $\alpha$ to $G_L$ becomes trivial. Hence the
central extension splits after going to  a finite extension, and so we
obtain a splitting $\rho_1': G_L\to H_1'(\bar{F})$ satisfying $f\circ
\rho_1'=\rho_1\!\mid G_L$. 

Since $\bar{\theta}\circ\bar{\rho}_1=\bar{\rho}_2$, we obtain  
\[ \pi\circ\theta\circ\rho_1'=(\pi\circ\rho_2)\!\mid G_L,\]
as morphisms from $G_L$ to $M(F)$. Since the kernel of $\pi$ at the
level of $F$-points is $Z_R(F)$, we obtain the first part of the theorem. 

To prove the  second and the third parts of the theorem, since
$Z_R\subset Z$, it follows from reasoning as above that there is an
isomorphism $\bar{\theta}$ between the images of $H_1$ and $H_2$ in
$M/Z$. Using the fact that $H_2/H_2\cap Z\simeq H_2Z/Z$, we now apply  
the second part of Lemma \ref{lift} to obain the proof in these
cases. 
\end{proof}

Part ii b) of the foregoing theorem allows us to avoid base
changing our original representation to an extension of the base
field, but instead allow the values of the twisting character to lie
in a connected, central group containing  $Z_R$.  
We now give a criterion for extending an isomorphism between two
representations.
\begin{lemma}\label{extension}
 Let  $\rho_1,  \rho_2:G_K\to GL_n(F)$ be $l$-adic
representations of $G_K$. Suppose there is a finite extension $L$ of
$K$ such that $\rho_1|G_L\simeq \rho_2|G_L$. Assume further that
$\rho_1|G_L$ is absolutely irreducible. Then there is a character
$\chi:G_K\to F^*$ such that,
\[\rho_2\simeq  \rho_1\otimes \chi.\]
\end{lemma}

\begin{proof} By Schur's lemma, the commuatant of $\rho_1({G_L})$ inside
the algebraic group $GL_n/F$,  is a form of $GL_1$.
By Hilbert Theorem 90, the commutant of $\rho_1(G_L)$ inside 
$GL_n(F)$, consists of precisely the scalar
matrices. 

For $\sigma\in G_K$, let $T(\sigma)=\rho_1(\sigma)^{-1}\rho_2(\sigma)$.
Since $\rho_1\!\mid_{G_L}=\rho_2\!\mid_{G_L}$, $T(\sigma)=Id$ for $\sigma\in
G_L$. Now for $\tau\in G_L$ and $\sigma\in G_K$,
\begin{equation*}
\begin{split}
T(\sigma)^{-1}\rho_1(\tau)T(\sigma) & =
\rho_2(\sigma)^{-1}\rho_1(\sigma)\rho_1(\tau)\rho_1(\sigma)^{-1}\rho_2(\sigma)
\\
&= \rho_2(\sigma)^{-1}\rho_1(\sigma\tau\sigma^{-1})\rho_2(\sigma) \\
& =\rho_2(\sigma)^{-1}\rho_2(\sigma\tau\sigma^{-1})\rho_2(\sigma)\\
& = \rho_2(\tau)\\
&=\rho_1(\tau)
\end{split}
\end{equation*}
Thus  $T(\sigma)$ is equivariant with respect
to the representation $\rho_1\!\mid_{G_L}$, and hence is given by a scalar
matrix $\chi(\sigma)$. Since $\chi(\sigma)$ is a scalar matrix, it follows
that for $\sigma,~\tau\in G_K$,
$\chi(\sigma\tau)=\chi(\sigma)\chi(\tau)$, i.e., $\chi$ is character
of $\mbox{Gal}(L/K)$ into the group of invertible elements $F^*$ of $F$,
and $\rho_2(\sigma)=\chi(\sigma)\rho_1(\sigma)$ for all $\sigma \in
G_K$.\\
\end{proof}

\begin{remark}
In the context of Question \ref{question}, the above lemma is
applicable,  whenever $H_{1d}$ has a
unique representation into $GL_n$, or as we shall see in the next
section when $R$ is a tensor or symmetric power of the original
representation and $L$ is such that the algebraic envelope of the
image of $(\rho_1\times \rho_2)|G_L$ is connected. In conjunction with
Theorem \ref{recoverthm}, it allows us to impose a mild hypothesis- that
$T$ has positive upper density for an answer to Question \ref{question}. 
\end{remark}

We now single out the case when $M=GL_n$. Other interesting examples
can be obtained by specialising $M$ to be $GSp_n$ or $GO_n$. For
the sake of simplicity we assume that $H^0_{1d}$ is simply connected.
 Observe that in the special
case when $\rho_1$ is further irreducible, by Remark \ref{rationality}
following Lemma \ref{lift}, 
 we obtain that $\theta$  is defined over $F$. 

\begin{corollary} Let $M=GL_n$, and assume  that the
semisimple component $H^0_{1d}$
of the connected component of the  algebraic envelope of the image of
$\rho_1$ is   simply connected. Let  $R:GL_n\to GL_m$ be a rational
representation, $c_i=|R(H_i)/R(H_i)^0|$, and let 
\[ T =\{v \in \Sigma_K-S\mid  {\rm Tr}({R\circ\rho_1}(\sigma_v))= {\rm Tr}({R\circ
\rho_2}(\sigma_v))\}.\]

a) Suppose that  $R\circ \rho_1$ and $R\circ \rho_2$ satisfy the
following {\bf DH2}:
\[ ud (T)  >{\rm min}(1-1/c_1,1-1/c_2). \]
Then there exists  an algebraic homomorphism $\theta:
H_1\to M$ defined over $F$,  a finite extension $L$ of $K$, and a
 character $\chi:
G_L\to GL_1(F)$ such that, 
\[ \rho_2\!\mid G_L=(\theta \circ \rho_1\!\mid G_L)\otimes \chi.\]

 b) Assume further that $H_1$ is connected and  
 that  $R\circ\rho_1$ and $R\circ\rho_2$ satisfy the following {\bf DH1}:
\[ ud (T) >1-1/2m^2. \]
Then there is an algebraic homomorphism $\theta:
H_1\to M$ defined over $F$, and a
 character $\chi:
G_K\to GL_1(F)$ such that, 
\[ \rho_2=(\theta \circ \rho_1)\otimes \chi.\]

c) If further  $H^0_{1d}$ has a unique absolutely irreducible representation
 to $GL_n$ upto equivalence,  and  if $ud(T)$ is positive, we have 
\[ \rho_2=\rho_1\otimes \chi.\]
\end{corollary}

\begin{remark} Let $G$ be a semisimple algebraic group over $F$, and
let $\theta_1,\theta_2,\cdots, \theta_m$ be the irreducible
representations of $G$ to $GL_n$.  Let $R$ denote the tensor product
$(F^n)^{\otimes m}$, on which $G$ acts via the tensor product of the
representations $\rho_1\otimes \cdots \otimes \rho_m$.  Suppose 
$\rho:G_K\to GL_n(F)$ is a Galois representation 
 with the Zariski closure of the image being
isomorphic to $G$. We do have for any $i$ that $R\circ \rho\simeq R\circ
(\theta_i\circ \rho)$. Thus if we are dealing with $R$ general, it is
not possible to limit $\theta$ further. In the next section, we will
see that if we restrict $R$ to be the tensor or symmetric power of the
original representation, then for `general' representations (meaning
non-induced),  we can
conclude under the hypothesis of the Theorem that $\theta$ can be
taken to be conjugation by an element of $GL_n(F)$. 
\end{remark} 

\begin{remark} In the automorphic context, if we replace $G_K$ by the
conjectural Langlands group $L_K$ whose admissible representations to
the dual $L$-group of a reductive group $G$ over $K$, 
parametrize  automorphic representations of $G(\A_K)$, 
 then corresponding to $R$ as
above, there is a conjectural functoriality lift. In this context the
theorems above and in the next section, amount to finding a
description of the fibres of the functoriality lifting with respect to
$R$. 
\end{remark}

\subsection{A finiteness result}
Given a representation $\rho: G_K\to M(F)$ and an element $m\in
M(\bar{F})$, define the twisted reprsentation $\rho^{(m)}:G_K\to
M(\bar{F})$ by 
\[\rho^{(m)}(\sigma)=m^{-1}\rho(\sigma)m.\]
Combining Theorem \ref{recoverthm} with a theorem of Richaradson
\cite{Ri}, we deduce  a finiteness
result upto twisting by representations with values in  the
center. For the sake of simplicity, we impose some additional
hypothesis on $\rho_1$.  

\begin{corollary} \label{finiteness}
Fix a  continuous representation 
$\rho_1:G_K\to M(F)$ as above. Assume that $H_1$ is connected,
$H_{1d}$ is simply connected, and that the center $Z$ of $M$ is
connected. Then there exists finitely many
representations $\rho_1, \rho_2, \cdots, \rho_d$ with the following
property: suppose $\rho:G_K\to M(F)$ is a continuous representation as
above such that $R\circ \rho_1$ and $R\circ \rho$ satisfy {\bf DH1}, where
$R:M\to GL_m$ is some rational representation of $M$ with kernel
contained in the center $Z$ of $M$. Then
there exists $m\in M(\bar{F})$ and a representation $\chi:G_K\to
Z(\bar{F})$ with values in $Z$, such that 
\[ \rho\simeq \rho_i^{(m)}\otimes \chi.\]
\end{corollary}

\begin{proof} By Theorem \ref{recoverthm}, we obtain that there exists
a morphism $\theta:H_{1d}\to M$ defined over $\bar{F}$, and a 
representation $\chi:G_K\to Z(\bar{F})$, such that under the
hypothesis of the Theorem, 
\[ \rho\simeq \rho_i^{(m)}\otimes \chi.\]
It follows from a theorem of Richardson \cite{Ri}, that upto conjugacy
by elements in $M(\bar{F})$, there exists only finitely many
representations of $H_{1d}$ into $M$. Choosing a representative in
each equivalence class, we obtain the Theorem. 
\end{proof}

\section{Tensor and Symmetric powers, Adjoint, Asai \\
 and Induced representations}
 We now specialize $R$ to some familiar
representations. The algebraic formulation of our problem, allows us
the advantage to  base change  and work over the complex numbers. 
For a linear representation $\rho$ of a group $G$ into
$GL_n$, let $T^k(\rho), ~S^k(\rho),
~E^k(\rho)~(k\leq n), ~{\rm Ad}(\rho)$ be respectively the $k^{th}$ tensor,
symmetric, exterior product and adjoint  representations of $G$. 

\subsection{Tensor powers}
\begin{proposition} \label{tensorconn}
Let $G$ be a connected algebraic group over a
characteristic zero base field $F$, and let $\rho_1, ~\rho_2$ be
finite dimensional semisimple representations of $G$ into $GL_n$. Suppose that 
\[ T^k(\rho_1)\simeq T^k(\rho_2)\]
for some $k\geq 1$. Then $\rho_1 \simeq \rho_2$. 
\end{proposition}

\begin{proof} We can work over $\C$. Let $\chi_{\rho_1}$ and
$\chi_{\rho_2}$ denote respectively the characters of $\rho_1$ and
$\rho_2$.  Since $\chi_{\rho_1}^k=\chi_{\rho_2}^k$, $\chi_{\rho_1}$
and $\chi_{\rho_2}$ differ by a $k^{th}$ root of unity. Choose a
connected neighbourhood $U$ of the identity in $G(\C)$, where the
characters are non-vanishing. Since $\chi_{\rho_1}(1)= \chi_{\rho_2}(1)$,
and the characters differ by a root of unity, we have $\chi_{\rho_1}=
\chi_{\rho_2}$ on $U$. Since they are rational functions on $G(\C)$
and $U$ is Zariski dense as $G$ is connected, we see that
$\chi_{\rho_1}=\chi_{\rho_2}$ on $G(\C)$. Since the representations
are semisimple, we obtain that $\rho_1$ and $\rho_2$ are equivalent.
\end{proof}

\begin{example} The connectedness assumption cannot be dropped. Fong
and Greiss \cite{fong-griess95} (see also Blasius\cite{blasius94}), 
have constructed for infinitely many triples $(n, q, m)$
homomorphisms  of 
$PSL_n(\F_q)$ into $PGL_m(\C)$,  which are elementwise conjugate but
not conjugate as representations. Here $\F_q$ is the finite field with
$q$-elements. Lift  two such homomorphisms  to representations
$$\rho_1, \rho_2 :SL(n, \F_q) \to GL(m,\C).$$
 We obtain that for
each $g\in SL(n,\F_q)$, $\rho_1(g)$ is conjugate to 
$\lambda \rho_2(g)$, with  $\lambda$ a scalar. Let $l$ be an exponent
 of the group $SL(n, \F_q)$. Then $g^l=1$ implies that
 $\lambda^l=1$. It follows that the characters $\chi_1$ and $\chi_2$
 of $\rho_1$ and $\rho_2$ respectively satisfy 
\[ \chi_1^l=\chi_2^l.\]
Thus the $l^{th}$ tensor powers of $\rho_1$ and $\rho_2$ are
equivalent, but by construction there does not exist a character
$\chi$ of $SL(n,\F_q)$ such that $\rho_2\simeq \rho_1\otimes \chi$.  
\end{example}

\subsection{Symmetric powers}
\begin{proposition} \label{symmetric}
Let $G$ be a connected reductive algebraic group
over a characteristic zero base field $F$. Let  $\rho_1, \rho_2$ be
finite dimensional representations of $G$ into $GL_n$. Suppose that 
\[ S^k(\rho_1)\simeq S^k(\rho_2)\]
for some $k\geq 1$. Then $\rho_1 \simeq \rho_2$. 
\end{proposition}

\begin{proof} We can work over $\C$. Let $T$ be a maximal torus of
$G$. Since two representations of a reductive group are equivalent if
and only if their collection of weights with respect to a maximal
torus $T$ are the same, it is enough to show that
 the collection of weights of $\rho_1$
and $\rho_2$ with respect to $T$ are the same. By Zariski density or
Weyl's unitary trick, we can work with a compact form of $G(\C)$ with
Lie algebra ${\mathfrak g}_r$. Let ${\mathfrak t}$ be a maximal torus inside
${\mathfrak g}_r$. The weights of the corresponding Lie algebra
representations associated to ${\rho_1}$ and $\rho_2$ are real
valued. Consequently we can order them with respect to a lexicographic
ordering on the dual of ${\mathfrak t}$. 

Let $\{\lambda_1,\cdots, \lambda_n\}$ (resp. $\{\mu_1,\cdots, \mu_n\}$) be
the weights of $\rho_1$ (resp. $\rho_2$) with $\lambda_1\geq
\lambda_2\geq \cdots\geq \lambda_n$ (resp. $\mu_1\geq \cdots \geq
\mu_n$). 
The weights of $S^k(\rho_1)$ are composed of elements of the form 
$$\left\{ \sum_{1\leq i\leq n} k_i\lambda_i \!\mid \sum_{1\leq i\leq n}
k_i=k\right\},$$ 
and similarly for $S^k(\rho_2)$.  

By assumption the weights of $S^k(\rho_1)$ and $S^k(\rho_2)$
are same. Since $k\lambda_1$ (resp. $k\mu_1$) is the highest weight of
$S^k(\rho_1)$ (resp. $S^k(\rho_2)$) with respect to the lexicographic
ordering, we have $k\lambda_1=k\mu_1$. Hence $\lambda_1=\mu_1$. 

By induction, assume that for $j<l, ~\lambda_j=\mu_j$. Then the set of
weights $\{\sum_{i<j} k_i\lambda_i\!\mid \sum k_i =k\} $ and 
$\{\sum_{i<j} k_i\mu_i\!\mid \sum k_i =k\} $ are same. Hence the
complementary  sets
$T_1(i)$ (resp. $T_2(i)$) composed of
weights in $S^k(\rho_1)$ (resp. $S^k(\rho_2)$), where at least one
$\lambda_j$ (resp. $\mu_j$) occurs with positive coefficient for some
$j\geq i$ are the same. 

The highest weight in $T_1(i)$ is
$(k-1)\lambda_1+\lambda_i$, and in $T_2(i)$ is $(k-1)\mu_1+\mu_i$. Since
$\lambda_1=\mu_1$, we obtain $\lambda_i=\mu_i$. Hence we have shown
that the collection of weights of $\rho_1$ and $\rho_2$ are the same,
and so the representations are equivalent. 
\end{proof}

\subsection{Adjoint and Generalized Asai representations}
 For a $G$-module $V$,
let ${\rm Ad}(V)$ denote the adjoint $G$-module given by the natural
action of $G$ on ${\rm End}(V)\simeq V^*\otimes V$,  
 where $V^*$ denotes the dual
of $V$. 

\begin{example} Let $G$ be a semisimple group, and let $V,W$ be non
self-dual irreducible representations of $G$.  Let $V_1=V\otimes W, ~V_2 =V\otimes
W^*$, considered as $G\times G$-modules.  Then as $G\times G$ irreducible modules (or as reducible $G$ modules), ${\rm Ad}(V_1)\simeq {\rm
Ad}(V_2)$, but $V_1$ is neither isomorphic to $V_2$ or to the dual
$V_2^*$. 
\end{example}

However when $G$ is simple, the following proposition is proved in
  \cite{Ra3}.
\begin{theorem}\label{tensorsimple}
Let $\fg$  be a simple Lie algebra over $\C$. Let $V_1, \cdots, V_n$
and $W_1, \cdots, W_m$ be non-trivial, irreducible, finite dimensional
$\fg$-modules. Assume that there is an isomorphism of the tensor
products, 
\[  V_1\otimes \cdots \otimes V_n\simeq W_1\otimes \cdots \otimes
W_m,\]
as $\fg$-modules. Then $m=n$, and there is a permutation $\tau$ of the
set $\{1, \cdots, n\}$, such that 
\[ V_i\simeq W_{\tau(i)},\]
as $\fg$-modules.

In particular, if $V, ~W$  are rreducible $\fg$-modules, and assume
that  
\[ {\rm End}(V)\simeq {\rm End}(W),\]
as $\fg$-modules. Then $V$ is either isomorphic to $W$ or the dual
$\fg$-module $W^*$. 
\end{theorem}

The above theorem arose out of the application to recovering
representations knowing that the adjoint representations are
isomorphic.  The arithmetical 
application of the above theorem 
to Asai representations was suggested by D. Ramakrishnan's
work, who has proved a similar result as below for the usual degree two
Asai representations, and I thank him for conveying to me his
results. 
We  consider now  a generalisation of Asai
representations.  Let $K/k$ be a Galois extension with Galois group
$G(K/k)$. Given $\rho$, we can associate the {\em pre-Asai representation}, 
\[As(\rho)=\otimes_{g\in G(K/K)}\rho^g,\]
where $\rho^g(\sigma)=\rho(\tilde{g}\sigma\tilde{g}^{-1}), ~\sigma \in
G_K$, and where $\tilde{g}\in G_k$ is a lift of $g\in G(K/k)$. At an
unramified place $v$ of $K$, which is split completely over a place
$u$ of  $k$, the
Asai character is given by,
\[ \chi_{As(\rho)}(\sigma_v)=\prod_{v|u}\chi_{\rho}(\sigma_v).\]
Hence we get that upto isomorphism, $As(\rho)$ does not
depend on the choice of the lifts $\tilde{g}$. If further $As(\rho)$
is irreducible, and $K/k$ ic cyclic, then $As(\rho)$ extends to a
representation of $G_k$ (called the Asai representation associated to
$\rho$ when $n=2$ and $K/k$ is quadratic).     

\subsection{Exterior powers}
 Let $V$ and $W$ be $G$-modules. It does not seem
possible to conclude in general from an equivalence  of the form
$E^k(V)\simeq E^k(W)$ as $G$-modules, that $V\simeq W$.  Here $E^k(V)$
denotes the exterior $k^{th}$ power of $V$. For example, let $V$ be a
non self-dual $G$-module of even dimension $2n$. If $G$ is semisimple,
then $E^n(V)$ is self-dual, and also that $E^n(V)$ is dual to
$E^n(V^*)$, where $V^*$ denotes the dual of $V$. Hence we have
$E^n(V)\simeq E^n(V^*)$, but $V\not\simeq V^*$. 

 It would be interesting to know, for the
possible applications in geometry,   the
relationship between two linear representations of  a connected
reductive group $G$ over $\C$ into
$GL_n$, given that their exterior $k^{th}$ power representations for some
positive integer $k<n$ are isomorphic.

\subsection{} We summarize the results obtained so far when $R$ is
taken to be a special representation.  

\begin{theorem}\label{tensorandsym} Let $\rho_1,~\rho_2:~G_K \to GL_n(F)$ be $l$-adic
representations of $G_K$.  Let $R$ be the representation 
$T^k$ or $S^k$  of $GL_n$, for some $k\geq 1$.
 Suppose that the representations $R\circ
\rho_1$ and $ R\circ\rho_2 $ satisfy {\bf DH2}. 
Then 

a) There is a finite extension $L$ of $K$ such that 
\[\rho_1\!\mid G_L \simeq \rho_2 \!\mid G_L.\]

b) Assume further that $H_1^0$ acts irreducibly on $F^n$. 
Then there is a character
$\chi:G_K\to F^*$ such that,
\[\rho_2\simeq  \rho_1\otimes \chi.\]

c) Suppose now that $H_1^0$ acts irreducibly on $F^n$ and 
 that the Lie algebra of the 
semisimple component of $H_1^0$ is simple. Assume now
 that the representations ${\rm Ad}(\rho_1)$ and  ${\rm Ad}(\rho_2)$
 satisfy {\bf DH2}. Then there is a character
$\chi:G_K\to F^*$ such that
\[\rho_2\simeq \rho_1\otimes \chi\quad \text{ or} \quad  
\rho_2\simeq \rho_1^*\otimes \chi.\]

d) Suppose now that $H_1^0$ acts irreducibly on $F^n$ and 
 that the Lie algebra of the  semisimple component of $H_1^0$ is simple.
Let $G(K/k)$ be a finite group of automorphisms acting on $K$, with
 quotient field $k$.  Assume now
 that the generalized Asai representations 
${\rm As}(\rho_1)$ and  ${\rm As}(\rho_2)$
 satisfy {\bf DH2}. Then there is a character
$\chi:G_K\to F^*$, and an element $g\in G(K/k)$ such that
\[\rho_2\simeq \rho_1^g\otimes \chi.\]

e) Suppose  that  the representations  $T^k(
\rho_1)$ and $T^k(\rho_2 )$ satisfy {\bf DH1}. Then if $\rho_1$ is
irreducible, so is $\rho_2$. 
\end{theorem}

\begin{proof} Let $\rho:=\rho_1\times \rho_2 :G_K\to (GL_n\times
GL_n)(F)$ be the product representation, and let
$H=H_{\rho}$. Applying the above Propositions \ref{tensorconn},
\ref{symmetric} and \ref{tensorsimple} to $H^0$ and Lemma
\ref{extension}, we obtain  parts a),
b), c) of the 
theorem. Part d) follows, since by assumptions the twisted tensor
product can be considered as different representations of the same
algebraic group (since the algebraic envelopes of the twists are
isomorphic), and applying Theorem \ref{tensorsimple}.   
Part e) follows from  Proposition \ref{tensorirr}  proved below. 

\end{proof}

\subsection{Induced representations}
Theorem \ref{tensorandsym} shows that it is easier to understand 
the generic situation when the algebraic
envelope is connected and the representation is absolutely
irreducible. The troublesome case then is to consider representations
which are induced from a subgroup of finite index. We now discuss
Question \ref{question} in this context, and hope that this will help
in clarifying some of the issues pertaining to endoscopy in the
context of the work of Labesse and Langlands and Blasius \cite{blasius94}.    
We first present another application of the algebraic machinery.
\begin{proposition}\label{tensorirr}
 Let $G$ be an algebraic group over a
characteristic zero base field $F$, and let $\rho_1, \rho_2$ be
finite dimensional semisimple representations of $G$ into $GL_n$. 
Suppose that \[ T^k(\rho_1)\simeq T^k(\rho_2)\]
for some $k\geq 1$. Then if  $\rho_1$ is irreducible, so is $\rho_2$. 
\end{proposition}
\begin{proof} The proof follows by base changing to complex numbers
and considering the characters $\chi_1,\chi_2$ 
of the representations restricted to a
maximal compact subgroup. The characters differ by a root of unity in
each connected component, and hence
$<\chi_1,\chi_1>=<\chi_2,\chi_2>$. Hence it follows by Schur
orthogonality relations that if one of the representations is irreducible,
so is the other.
\end{proof}

We continue with the hypothesis of the above proposition.   By Proposition
\ref{tensorconn}, we can assume  
that $\rho_1|G^0$ and  $\rho_2|G^0$ are isomorphic, which we will
denote by $\rho^{0}$ .  Write 
\[\rho^0=\oplus_{i \in I} r_i,\] where $r_i$ are the
representations on the isotypical components, and $I$ is  the
indexing set of the isotypical components of $\rho^{(0)}$. The group
$\Phi:=G/G^0$ acts via the representations $\rho_1$ and $\rho_2$ to give
raise to two permutation representations $\sigma_1$ and $\sigma_2$ on
$I$. For each $\phi\in \Phi, ~l=1,2$ let 
\[ S_l(\phi)=\{ i\in I\mid \sigma_l(\phi)(i)=i\}.\]
Denote by $V_l^{\phi}$ the representation space of $H^0$ obtained by
taking the direct sum of the representations indexed by the elements 
occuring in $S_l(\phi)$. 
For $\phi\in \Phi$, let $G^{\phi}$ denote the corresponding connected
component (identity component is $G^0$). 
Let $p_l^{\phi}, ~p_l^0$ denote respectively the projections of
$G^{\phi},~ G^0$   to $GL(V_l^{\phi})$. 
With these assumptions we have, 
\begin{proposition}\label{tensorinduced} 
 Let $G$ be an algebraic group over a
characteristic zero base field $F$, and let $\rho_1, \rho_2$ be
finite dimensional semisimple representations of $G$ into $GL_n$. 
Suppose that \[ T^k(\rho_1)\simeq T^k(\rho_2)\]
for some $k\geq 1$. Assume further that either $\rho_1$ or $\rho_2$ is
absolutely irreducible. 
With notation as above, we have for all $\phi\in \Phi$, $S_1(\phi)=S_2(\phi)$. In
particular we have ${\rm Ker}(\sigma_1)={\rm Ker}(\sigma_2)$. 

Assume further that one of the representations is irreducible. Then
$\rho_1$ and $\rho_2$ are induced respectively from representations
$r_1', ~r_2'$ of  the same subgroup $G'$ of $G$, such that the
restriction of $r_1'$ and $r_2'$ to $G^0$ is an isotypical component
of $\rho^0$.   
\end{proposition} \label{inducedsubgroup}
\begin{proof}  Let $\eta_i$ denote the
character of $r_i$. Let $G^{(\phi)}$ denote the coset of $G^0$ in $G$
corresponding to $\phi\in \Phi$. Denote by $\chi_1$ and $\chi_2$ the
characters of $\rho_1$ and $\rho_2$ respectively.  For an element
$g\phi\in G^{\phi}, \, g\in G^0, \, l=1,2$ we have 
\[ \chi_l(g\phi)=\sum_{i\in S_l(\phi)} \eta_i(g).\]
Since $\chi_1^k=\chi_2^k$, we obtain 
\[\sum_{i\in S_1(\phi)} \eta_i=\zeta\sum_{i\in S_2(\phi)} \eta_i,\]
for some $\zeta$ a $k^{th}$ root of unity. The first part of the 
lemmma follows from 
linear independence of irreducible characters of a group. 

To prove the second assertion, fix an isotypical component say $r_{i_0}$ of
$\rho^0$. The subgroup $G'$ of $G$ which stablizes $r_{i_0}$ is the
same for both the representations. Let $r_1', ~r_2'$ be the extension
of $r_{i_0}$ as 
representations of $G'$ associated respectively to $\rho_1$ and
$\rho_2$. It follows by theorems of Clifford that
the representations $\rho_1$ and $\rho_2$ are induced from $r_1'$ and
$r_2'$ respectively.
\end{proof}

Assume now that the constituents $r_i$ of
$\rho^0$ are absolutely irreducible, i.e, the irreducible
representations occur with multiplicity one.
In the notation of Proposition \ref{tensorinduced}, the assumption of
multiplicity one on $\rho^0$, implies that $r_2'=r_1'\otimes \chi$ for
some character $\chi\in {\rm Hom}(G', F^{'*})$ trivial upon restriction
to $G^0$, where $F'$ is a finite extension of $F$. 
 Assuming the hypothesis of Proposition \ref{tensorirr}, we
would like to know whether $\rho_2$ and $\rho_1$ differ by a
character. This amounts to knowing that the character $\chi$ extends
to a character of $G$, since the representations are induced. 
Assume from now onwards that $G'$ is normal in $G$. Then the question
of extending $\chi$ amounts first to showing that $\chi$ is invariant
and then to show that invariant characters extend.

\begin{remark}
We rephrase the problem in a different language, with the hope that it
may shed further light on the question. For $\sigma\in G$, let
$T(\sigma)=\rho_1(\sigma)^{-1}\rho_2(\sigma)$ be as in Lemma
\ref{extension}. The calculations of Lemma \ref{extension} with now
$\sigma \in G, \tau \in G^0$, show that $T(\sigma)$ takes values in
the commutant of $\rho^0$. Let $S$ denote the
$F'$-valued points of 
commutant torus of $\rho^0$ which is defined over a finite extension  $F'$ of
$F$.  
Since we have assumed that $\rho_1$ is
irreducible, $G$ acts on $S$ via $\sigma_1$ transitively as a
permutation representation  on 
indexing set $I$. Hence $I$ can be taken to be $G/G'$ and 
$S$ is isomorphic to the induced module ${\rm Ind}_{G'}^G(F^{'*})$, where
the action of $G'$ on $F^{'*}$ is trivial. We have for $\sigma, \tau \in
G$, 
\[ T(\sigma\tau)=\rho_1(\tau)^{-1}T(\sigma)\rho_1(\tau)T(\tau), \]
i.e, $T$ is a one cocycle on $G$ with values in $S$. 
Since $S$ is induced, the `restriction' map $H^1(G, S)\to H^1(G', F^{'*})$
is an isomorphism. The invariants of the $G$-action on $S$ is given by
the diagonal $F^{'*}$ sitting inside $S$, and the composite map,
$ {\rm Hom}(G, F^{'*}) \to H^1(G, S)\to {\rm Hom}(G', F^{'*})$ 
is the restriction
map. To say that $\chi$ extends to a character of $G$, amounts to
knowing that $\chi$ lies in the image of this composite map. 
\end{remark}

\begin{example} We consider the example given by Blasius \cite{blasius94} in
this context. Let $n$ be an odd prime, and let $H_n$ be the  finite
Heisenberg group, with generators $A, ~B, ~C$ subject to  the
relations: $A^n=B^n=C^n=1, ~AC=CA, ~BC=CB, ~AB=CBA$. Let $e_1, \cdots,
e_n$ be a basis for $\C^n$, and let $\xi_n$ be a primitive $n^{th}$
root of unity. For each integer
$a$ coprime to $n$, define the representation $\rho_a:H_n\to GL_n(\C)$
by, 
\begin{equation*}
\begin{split}
 \rho_a(A)e_i& =\xi_n^{(i-1)a}e_i\\
 \rho_a(B)e_i& =e_{i+1}\\
 \rho_a(C)e_i& =\xi_n^{a}e_i,
\end{split}
\end{equation*}
where the notation is that $e_{n+1}=e_1$. It can be seen that $\rho_a$
are irreducible representations, and that the corresponding projective
representations for any pair of integers $a, ~b$ not congruent modulo
$n$,  are inequivalent. Further for any element $h\in H_n$,
the images of $\rho_a(h)$ and $\rho_b(h)$  in $PGL(n, \C)$ are
conjugate. Hence it follows that for some positive integer $k$ (which
we can take to be $n$) the representations  $T^k(\rho_a)$ and
$T^k(\rho_b)$ are isomorphic. 

Let $T$ be the abelian normal
subgroup of index $n$ generated by $A$ and $C$. There exists a
character $\chi_{ab}$ of $T$ such that $\rho_a|T\simeq \rho_b|T\otimes
\chi_{ab}$. From the theory of
induced representations, it can be further checked that $\rho_a|T$ has
multiplicity one, and that $\rho_a$ is induced from a character
$\psi_a$ of $T$. However we have that there does not exist any
character $\eta$ of $H_n$ such that $\rho_a\simeq \rho_b\otimes \eta$.
\end{example}

With this example in mind, we now present proposition in the positive
direction (which applies in particular to CM forms of weight $\geq 2$): 

\begin{proposition}\label{inducedinv}
 Let $G$ be an algebraic group over a
characteristic zero base field $F$, and let $\rho_1, \rho_2$ be
finite dimensional semisimple representations of $G$ into $GL_n$. 
Suppose that \[ T^k(\rho_1)\simeq T^k(\rho_2)\]
for some $k\geq 1$. Assume further that either $\rho_1$ or $\rho_2$ is
absolutely irreducible, and the following assumptions: 
\begin{itemize}
\item the representation $ \rho^{0}:=\rho|G^0$ can be written as a
direct sum of irreducible representations $\oplus_{i\in I}r_i$  
with  multiplicity one, i.e.,
each of the isotypical components $r_i$ are irreducible. 

\item the subgroup $G'$ is normal in $G$. 
\end{itemize} 
Let $r_1', ~r_2'$ be representations of $G'$ as in the proof of
Proposition \ref{inducedsubgroup}. Let $\chi$ be a character of $G'$
such that $r_2'\simeq r_1'\otimes \chi$ (this exists because of the
assumption of multiplicity one). 

Then $\chi$ is invariant with respect to the action of $G$ on the
characters of $G'$. In particular, if invariant characters of $G'$
extend to invariant characters of $G$ (which happens if $G/G'$ is
cyclic), then we have 
\[ \rho_2\simeq \rho_1\otimes \chi.\]
\end{proposition}
\begin{proof}
Restricting $\rho_1$ and $\rho_2$ to $G'$, have by our assumptions
\[ \rho_1|G'=\oplus_{\substack{\phi \in G/G'}}{r_1'}^{\phi},\quad {\rm
and}\quad \rho_2|G'=\oplus_{\substack{\phi \in G/G'}}({r_1'\chi})^{\phi},\]
where $\chi$ is a character of $G'$ trivial on $G^0$.   Let $\mu_1$ denote the character of the representation $r_1$ of $G^0$.
Let $\tau$ be an element of $G-G'$. 
If $\chi\neq \chi^{\tau}$, choose
an element $\theta$ of $G'$ such that $\chi(\theta)\neq
\chi^{\tau}(\theta)$. It follows from our hypothesis of multiplicity
one,  that  any element $x$ in the coset $G^{\theta}$ of $G^0$ 
 defined by $\theta$ can be written as
$x=zy$, where $z$ is a fixed element  in the center of $G'$ depending
only on $\theta$ and not on $x$,  and $y \in G^0$. 
We obtain from our assumption 
$T^k(\rho_1)\simeq T^k(\rho_2)$,  that for some $\zeta \in {\mathbf
\mu}_k$, we have  
\[\zeta\sum_{\phi\in G/G'}(\mu_1(z)\mu_1(y))^{\phi}=\sum_{\phi\in
G/G'}(\mu_1(z)\mu_1(y)\chi(\theta))^{\phi}, \]
for all $y\in G^0$. But the assumption of multiplicity one implies
that the irreducible characters $\mu_1^{\phi}$ of $G^0$ are linearly
independent, and hence the above equality forces the character $\chi$
to be invariant.

\end{proof}

\section{Applications to Modular Forms and Abelian Varieties}
We consider now applications of the general theory developed so far to
the $l$-adic representations attached to holomorphic (Hilbert) modular
forms and to abelian varieties. As a corollary to Theorem
\ref{recoverthm} and Proposition \ref{inducedinv},  we have

\begin{corollary}  Let $K$ be a number field and let $\rho_1,
~\rho_2:G_K\to GL_2(F)$ be continuous $l$-adic representations as
above. Suppose that the algebraic envelope  $H_1$ of the image  $\rho_1(G_K)$
contains a maximal torus.  Let $R:GL_2\to GL_m$ be a rational representation
with kernel contained in the center of $GL(2)$, and $T$ be as in
Question \ref{question}.  
Then the following holds:

a)  Assume that the algebraic envelope of the image of $\rho_1$
contains $SL_2$. If $ud(T)$ is positive,
 then there exists a character 
$\chi:G_K\to GL_1(\bar{F})$ such that $\rho_2\simeq \rho_1\otimes
\chi$.

b)  Suppose that   $ ud(T) > 1/2$.   
Then there exists a finite extension $L$ of $K$ and a character $\chi:G_L\to GL_1(\bar{F})$
such that $\rho_2|G_L\simeq \rho_1|G_L\otimes \chi$. 

c) Suppose  that  the representations  $T^k(
\rho_1)$ and $T^k(\rho_2 )$ satisfy {\bf DH1}, for some positive
integer $k$ and  $\rho_1$ is
irreducible. Then there exists a character $\chi:G_K\to GL_1(\bar{F})$
such that $\rho_2\simeq \rho_1\otimes \chi$. 

\end{corollary}

We now make explicit the application to classical holomorphic modular
forms. A similar statement as below can be made for the class of
Hilbert modular forms too. For any  pair of positive integers $N,
k\geq 2$, and Nebentypus character $\omega: (\Z/N\Z)^*\to \C$, denote
by $S_k^n(N,\omega)$ the set of normalized newforms of weight $k$ on
$\Gamma_1(N)$. The assumption on the weights imply that the Zariski
closure of the image of the $l$-adic representation $\rho_f$ of
$G_{\Q}$  attached to $f$
contains a maximal torus, thus 
satisfying  the hypothesis of the above corollary. For $f\in  S_k^n(N,\omega)$ and $p$ coprime to $N$, let
$a_p(f)$ denote the corresponding Hecke eigenvalue. 

\begin{corollary} Let $f\in S_k^n(N,\omega)$ and $f'\in
S_k^n(N',\omega')$ and let $\rho_f,\rho_{f'}$ be the associated
$l$-adic representations of $G_{\Q}$.  Let $R:GL_2\to GL_m$ be a rational representation
with kernel contained in the center of $GL(2)$, and $T$ be as in
Question \ref{question}.  
Then the following holds:

a) If $f$ is not a CM-form and  $ud(T)$ is positive,
 then there exists a Dirichlet  character 
$\chi$ such that for all $p$ coprime to $NN'$, we have 
\[ a_p(f)=a_p(f')\chi(p).\]
In particular $k=k'$ and $\omega=\omega'\chi^2$. 

b) Suppose  that for some positive integer $l$, we have 
\[  a_p(f)^l=a_p(f')^l,\]
for a set of primes of density strictly greater than
 $1-2^{-(2l+1)}$. Assume also that if $f$ is dihedral, then it is not
 of weight one. 
 Then there exists a Dirichlet character $\chi$
such that for all $p$ coprime to $NN'$, we have 
\[ a_p(f)=a_p(f')\chi(p).\]

If we further assume that the conductors $N$ and $N'$ are squarefree
in a) and b) above, we can conclude that $f=f'$. 

\end{corollary}

This gives us the generalization to classical modular forms of the
theorems of Ramakrishnan \cite{DR}. We now consider an application to
Asai representations associated to holomorphic Hilbert modular
forms. With notation as in Part d) of Theorem \ref{tensorandsym}, assume
further that $K/k$ is a quadratic extension of totally real number
fields, and that $\rho_1, ~\rho_2$ are $l$-adic representations
attached  respectively to holomorphic Hilbert modular forms $f_1, ~f_2$
over $K$. Let $\sigma$ denote a generator of the Galois group of
$K/k$, and further assume that $f_1$ (resp. $f_2$)  is not isomorphic to any twist of
$f_1^{\sigma}$ (resp. $f_2^{\sigma}$ by a character, where
$f_1^{\sigma}$ denotes the form $f_1$ twisted by $\sigma$. Then the
twisted tensor  automorphic representations of $GL_2(\A_K)$, defined
by $f_i\otimes f_2^{\sigma}$ (with an abuse of notation) is
irreducible. Further it has been shown to be modular by
D. Ramakrishnan \cite{DR2}, and descend to define automorphic
representations denoted respectively by  
 $As(f_1)$ and $As(f_2)$ on $GL_2(\A_k)$. As a simple consequence of
Part d) of Theorem \ref{tensorandsym}, we have the following corollary:

\begin{corollary} Assume further that $f_1$ and $f_2$ are non-CM
forms, and that the Fourier coefficients of $As(f_1)$ and $As(f_2)$ are
equal at a positive density of places of $k$, which split in $K$. Then
there is an idele-class character $\chi$ of $K$ satisfying
$\chi\chi^{\sigma}=1$ such that,
\[ f_2=f_1\otimes \chi \quad \text{or}\quad  f_2 =f_1^{\sigma}\otimes
\chi.\]
\end{corollary}

We limit ourselves now to
providing a simple  application
of the theorems in the context of abelian varieties. Let
$K$ be a number field, and $A$ be  an abelian variety  of dimension $d$
over $K$,  and
let $\rho_A$ be the associated $l$-adic representation of
$G_K$ on the Tate module $V_l(A)={\varprojlim}        
A_{l^n}\otimes_{\Z}\Q\simeq GL_{2d}(\Q_l)$, where $A_{l^n}$ denotes
the group of $l^n$ torsion points of $A$.  
 As applications of our general theorems, we have
\begin{corollary} Let $A, B$ be abelian varieties as above. 

 a) Let $R$ be either $T^k$ or $S^k$  for some
positive integer $k$. Suppose that $R\circ \rho_A$ and $ R\circ
\rho_B$ satisfy the density hypothesis {\bf DH2}. Then $A$ and $B$ are
isogenous over $\bar{K}$. 

If further the connected component of the algebraic envelope $H_A$ of
$\rho_A$ acts absolutely irreducibly on the Tate module $V_l(A)$, then
there exists a Dirichlet character $\chi$ of $K$ such that
$\rho_B\simeq \rho_A\otimes \chi$. 

b) Suppose that we are in the `generic' case, i.e., $H_A =GSp_d$. Let 
$R$ be  any representation of $GSp_d$ with finite kernel, 
such  that ${\rm Tr}(R\circ \rho_A)(\sigma_v)={\rm Tr}(R\circ
\rho_B)(\sigma_v)$ on a set of places $v$ of $K$ of positive upper
density. Then  $A$ and $B$ are
isogenous over $\bar{K}$, and there exists a Dirichlet character $\chi$ of $K$ such that
$\rho_B\simeq \rho_A\otimes \chi$.
\end{corollary}

\end{document}